\documentclass[11pt]{amsart}
\usepackage[T1]{fontenc}
\usepackage{lmodern}
\usepackage{amsmath,amssymb,mathtools}
\usepackage[margin=1.15in]{geometry}
\usepackage{microtype}
\usepackage[colorlinks=true,linkcolor=blue,citecolor=blue,urlcolor=blue]{hyperref}
\hypersetup{pdftitle={Translation continuity for the local bilinear circular maximal operator},
 pdfauthor={Xuyi Wei},pdfsubject={Translation continuity and bilinear circular maximal averages},
 pdfkeywords={bilinear maximal operator, cinematic curvature, local smoothing}}
\numberwithin{equation}{section}
\newtheorem{theorem}{Theorem}[section]
\newtheorem{proposition}[theorem]{Proposition}
\newtheorem{lemma}[theorem]{Lemma}
\newtheorem{corollary}[theorem]{Corollary}
\theoremstyle{remark}
\newtheorem{remark}[theorem]{Remark}
\newcommand{\R}{\mathbb R}
\newcommand{\C}{\mathbb C}
\newcommand{\A}{\mathcal A}
\newcommand{\M}{\mathcal M_{\mathrm{loc}}}
\newcommand{\supp}{\operatorname{supp}}

\newcommand{\norm}[1]{\left\|#1\right\|}

\title[Continuity of bilinear circular maximal averages]{Translation continuity for the local bilinear circular maximal operator}
\author{Xuyi Wei}
\date{}

\begin{document}
\begin{abstract}
We prove a positive power translation-continuity estimate for the local
bilinear circular maximal operator on the real line. More precisely,
there is $\eta>0$ such that the operator applied to a translation
difference in either input maps $L^8(\R)\times L^8(\R)$ to $L^4(\R)$
with norm $O(|h|^\eta)$. This answers Question~6.3 in Borges's survey
on bilinear spherical averages. The main step is a frequency-localized
estimate for arbitrary bounded inputs. A $TT^*$ argument converts it
into a linear maximal problem on the plane, associated with pairs of
circles. We verify the corrected multi-parameter cinematic curvature
condition of Chen--Guo--Yang by proving an irreducibility and
infinitesimal independence statement for the resulting quartic curves.
Uniform algebraic sublevel estimates control the degeneracies.
Interpolation also gives continuity throughout the finite H\"older
range $2<p,q<\infty$, with output exponent $pq/(p+q)$.
\end{abstract}
\maketitle
\enlargethispage{3pt}

\section{Introduction}

Let $d\sigma$ denote normalized arclength measure on $S^1$. For $t>0$,
define
\begin{equation}\label{eq:average}
 \A_t(f,g)(x)=\int_{S^1}f(x-ty_1)g(x-ty_2)\,d\sigma(y)
 =\frac1{2\pi}\int_0^{2\pi}f(x-t\cos\theta)
                  g(x-t\sin\theta)\,d\theta.
\end{equation}
The local maximal operator is
\begin{equation}\label{eq:maximal}
 \M(f,g)(x)=\sup_{1\leq t\leq2}|\A_t(f,g)(x)|.
\end{equation}
Thus the supremum is taken \emph{inside} the output norm. We use the
translation convention $\tau_hf(x)=f(x-h)$.

In \cite[Question~6.3]{Borges}, Borges asks whether translation
continuity estimates for \eqref{eq:maximal} hold in dimension one.
The fixed-radius problem and the local maximal problem are different:
a bound for $\sup_t\|\A_t((I-\tau_h)f,g)\|_r$ does not answer this
question. We establish the following local maximal estimate.

\begin{theorem}\label{thm:main}
There are constants $\eta>0$ and $C<\infty$ such that, for every
$f,g\in L^8(\R)$ and every $h\in\R$,
\begin{equation}\label{eq:main}
 \norm{\M((I-\tau_h)f,g)}_4
 +\norm{\M(f,(I-\tau_h)g)}_4
 \leq C |h|^\eta\norm f_8\norm g_8.
\end{equation}
\end{theorem}

The exponent is positive but is not optimized. In particular, the
theorem asserts the existence requested in Question~6.3, rather than
an optimal continuity exponent or a complete range of Lebesgue indices.
The inputs in \eqref{eq:main} are arbitrary; no restriction to a
frequency configuration, a selected oscillatory phase, or a regular
radius selector is imposed.

The role of the triple $(8,8,4)$ is to provide one point at which a
positive continuity power is available. It is not a proposed sharp
range. Proposition~\ref{prop:interpolation} below gives the precise
interpolation principle: continuity at this point propagates to the
interior of the strong boundedness region. In particular, the known
H\"older bounds in dimension one give the following consequence.

\begin{corollary}\label{cor:holdercontinuity}
For every $2<p,q<\infty$, put $r=pq/(p+q)$. There are
$\eta_{p,q}>0$ and $C_{p,q}<\infty$ such that
\begin{equation}\label{eq:holdercontinuity}
 \norm{\M((I-\tau_h)f,g)}_r+
 \norm{\M(f,(I-\tau_h)g)}_r
 \leq C_{p,q}|h|^{\eta_{p,q}}\norm f_p\norm g_q
 \quad(0<|h|\leq1).
\end{equation}
\end{corollary}

The boundedness input for this corollary is the work of Christ--Zhou
\cite{CZ} and Dosidis--Ramos \cite{DR}. Its short deduction is given
after Proposition~\ref{prop:interpolation}. We do not assert continuity
at the boundary $p=2$ or $q=2$, nor a uniform power as either exponent
approaches a boundary.

We prove a stronger intermediate estimate. Fix a smooth multiplier
$p$ supported in $\{\xi:c\leq|\xi|\leq C_0\}$, where $c>0$, and write
$P_N=p(D/N)$ for $N\geq2$. The constants below may depend on this fixed
multiplier.

\begin{proposition}\label{prop:frequency}
There is $s>0$ such that
\begin{equation}\label{eq:globalfrequency}
 \norm{\M(P_Nf,g)}_4
 \leq C N^{-s}\norm f_8\norm g_8
 \qquad(N\geq2).
\end{equation}
The symmetric estimate with the projection in the second input also
holds. In particular, the other input need not be frequency localized.
\end{proposition}

The proof first establishes a local $L^\infty\times L^\infty\to L^1$
estimate with a positive power of $N^{-1}$. After linearizing the
radius supremum, Cauchy--Schwarz removes the second input. The kernel
that remains is a linear average of
$H(u,w)=F(u)\overline{F(w)}$ over a curve obtained by composing two
circle-coordinate maps. Its output has two independent coordinates.
This permits application of the variable-coefficient linear maximal
smoothing theorem of Chen--Guo--Yang \cite{CGY}.

The essential geometric point is that the four parameter derivatives
of a pair-circle graph are linearly independent on every regular
branch whenever the two centers differ. This follows from the
irreducibility of an explicit quartic polynomial. The resulting
cinematic determinant is therefore not identically zero for any
admissible pair of radii. An algebraic norm and a polynomial sublevel
estimate make this assertion quantitative, uniformly in the radii.

\subsection*{Related work}
The linear theory of circular maximal operators provides the analytic
background for the argument. Bourgain's work on planar curved averages
\cite{Bourgain} was followed by Sogge's variable-coefficient formulation
\cite{Sogge} and the local smoothing theory of Mockenhaupt, Seeger, and
Sogge \cite{MSS}. The multi-parameter theorem of Chen, Guo, and Yang
\cite{CGY} extends this cinematic-curvature framework to families with
several maximal parameters. Its proof combines variable-coefficient
decoupling, as developed by Beltran, Hickman, and Sogge \cite{BHS},
the bootstrapping method of Pramanik and Seeger \cite{PS}, and
moment-curve decoupling of Bourgain, Demeter, and Guth \cite{BDG}.
In this paper we use the corrected variable-coefficient hypothesis
in version~7 of \cite{CGY}, and record the quantitative dependence
needed for frequency-dependent excision in Appendix~\ref{app:quantitative}.

In the bilinear setting, Borges, Foster, Ou, Pipher, and Zhou
\cite{BFOPZ} developed continuity estimates as an ingredient in sparse
domination. Their one-dimensional fixed-radius continuity estimate
uses the trilinear smoothing inequality of Christ and Zhou \cite{CZ};
their localized maximal continuity estimates apply in dimensions at
least two. For the full bilinear circular maximal operator on the
line, Christ and Zhou \cite[Theorem~1.3(ii)]{CZ} and, independently,
Dosidis and Ramos \cite[Theorem~1]{DR} established the H\"older bounds
$L^p\times L^q\to L^{pq/(p+q)}$ for $p,q>2$. These boundedness results
provide interpolation endpoints, while the present argument supplies
a positive frequency gain with a continuum of radii. Borges's
survey \cite{Borges} places this continuity problem in the broader
study of bilinear spherical averages.

The algebraic part of our proof has a different role from the linear
analytic input. Irreducibility of the pair-circle quartic yields
independence of its parameter derivatives. The analytic Wronskian
criterion, recalled for example in \cite{BostanDumas}, converts this
independence into a nontrivial cinematic determinant. A compact
semialgebraic \L{}ojasiewicz inequality \cite{Coste} then makes its
nonvanishing quantitative, uniformly in the two radii. Finally,
multilinear interpolation in the form of
\cite[Proposition~2.2]{BOS} explains why one continuity triple suffices
to obtain continuity throughout the interior of a known boundedness
region.

\subsection*{Conventions}
We use the Fourier transform
$\widehat f(\xi)=\int_{\R}f(x)e^{-ix\xi}\,dx$ and the corresponding
inverse transform with factor $(2\pi)^{-1}$; in two dimensions the
inverse factor is $(2\pi)^{-2}$. A frequency annulus always has fixed
inner and outer ratios. Lebesgue norms without a domain are taken on
$\R$, except for functions explicitly defined on $\R^2$.
All unspecified compact coordinate intervals have fixed size.
Constants denoted by $C$ may increase from one occurrence to the next.
An assertion that a quantity is bounded polynomially in inverse
margins means a bound by a fixed power of their product. All such
powers are independent of $N$ and of the selected radius function.
We first work with smooth inputs and finite radius suprema; the
extension to the stated function spaces is justified in Section~\ref{sec:completion}.

\section{A quantitative linear smoothing input}\label{sec:linear}

We state the form of the linear result needed below, including its
dependence on nondegeneracy margins. For a graph family
$w=\phi(u;z_1,z_2,v_1,v_2)$, let
\begin{equation}\label{eq:jet}
 \mathfrak W(u;z,v)
 =\det\left[\partial_u^j\partial_{p}\phi(u;z,v)
       \right]_{\substack{0\leq j\leq3\\p\in(z_1,z_2,v_1,v_2)}}.
\end{equation}
This is the determinant of the four parameter derivatives of the
graph jet of orders zero through three.

\begin{lemma}[Quantitative cinematic maximal smoothing]\label{lem:linear}
Fix $p>12$. There are an integer $m$, a constant $C_p$, and a finite
exponent $K$ with the following property. Consider
\begin{equation}\label{eq:linearoperator}
 \mathcal S_v H(z)=\int H(\Gamma(z,t;v))a(z,t;v)\,dt,
 \qquad z,v\in\R^2,
\end{equation}
on fixed compact coordinate boxes. Suppose that the derivatives
through order $m$ of the curve and amplitude, the $C^m$ norms of the
regular graph-coordinate changes and their inverses, and the
reciprocals of the relevant support margins are bounded by
$A\geq2$. Suppose also that
\begin{equation}\label{eq:H12}
 |\det \partial_z\Gamma|\geq A^{-1},\qquad
 |\partial_t\Gamma|\geq A^{-1},\qquad
 |\mathfrak W|\geq A^{-1}
\end{equation}
on the amplitude support, with fixed powers of $A$ allowed in these
assumptions. If $\widehat H$ is supported in a fixed annulus of radius
$N$, then
\begin{equation}\label{eq:lineardecay}
 \norm{\sup_v|\mathcal S_vH|}_{L^p_z}
 \leq C_p A^K N^{-1/(2p)}\norm H_{L^p(\R^2)}.
\end{equation}
A polynomial number of regular coordinate patches may be summed by
increasing $K$.
\end{lemma}

This lemma is the finite-jet quantitative form of
\cite[Theorem~1.5 and (1.39)]{CGY}. We use \emph{version~7}, with its
corrected curvature hypothesis. For clarity, in the normalized
coordinates of that paper,
\[
 \Gamma(X,Y;\theta;v)=(X-\theta,h(X,Y;\theta;v)),\qquad
 T=(h,h_\theta,h_{\theta\theta},h_{\theta\theta\theta})^T,
\]
the corrected determinant is
\begin{equation}\label{eq:correctedH3}
 \det[T_\theta+T_X,T_Y,T_{v_1},T_{v_2}].
\end{equation}
It includes the zeroth graph jet. Set $u=X-\theta$. The vector field
$\partial_\theta+\partial_X$ differentiates the parameters while
holding $u$ fixed. Consequently \eqref{eq:correctedH3} is precisely
\eqref{eq:jet}, up to the nonzero factors of the coordinate changes
and the row signs $(-1)^j$. This observation will be used repeatedly.

The polynomial factor in \eqref{eq:lineardecay} is needed because the
good regions below depend on $N$. Appendix~\ref{app:quantitative}
records the normalization argument giving this factor from the proof
in \cite{CGY}. In particular, we do not infer polynomial dependence
from a qualitative compactness argument.

\section{The geometry of two circles}\label{sec:geometry}

\subsection{The composition quartic}
Let $a\ne b$ be real and let $r,s>0$. Consider pairs $(u,w)$ for which
there is a common intermediate coordinate $v$ satisfying
\begin{equation}\label{eq:twocircles}
 (u-b)^2+(v-b)^2=s^2,\qquad
 (w-a)^2+(v-a)^2=r^2.
\end{equation}
Put $z=u-b$, $y=w-a$, $h=b-a$, $R=r^2$, $S=s^2$, and
$c=R-S-h^2$. Subtracting the equations in \eqref{eq:twocircles} gives
$2h(v-b)=c+z^2-y^2$. Hence the pair curve is exactly
\begin{equation}\label{eq:quartic}
 F(u,w;a,b,R,S)
 =(z^2-y^2+c)^2-4h^2(S-z^2)=0.
\end{equation}
Conversely, a solution of \eqref{eq:quartic} determines $v$ by the
displayed linear equation and satisfies both circle equations.

\begin{lemma}\label{lem:irreducible}
If $a\ne b$ and $R,S>0$, the polynomial in \eqref{eq:quartic} is
irreducible over $\C$.
\end{lemma}
\begin{proof}
Over $K=\C(z)$, adjoin $q=\sqrt{S-z^2}$. This extension has degree
two, since $S-z^2$ has two simple zeros. The other square root obeys
$y^2=z^2+c-2hq$. The norm of its radicand down to $K$ is
\begin{equation}\label{eq:normquartic}
 z^4+(2c+4h^2)z^2+c^2-4h^2S.
\end{equation}
This is not a square in $K$. Indeed, a polynomial that is a square
in $K$ is a square in $\C[z]$: if a rational square has no finite
poles, its square root has no finite poles and hence is a polynomial.
Comparing the leading, cubic, and
quadratic coefficients would force the square to be
$(z^2+c+2h^2)^2$. Equality of the constant coefficients would then
give $4h^2(c+h^2+S)=4h^2R=0$, a contradiction.
Thus $y^2$ is not a square in $K(q)$, since the norm of a square is a
square. Moreover $q=(z^2+c-y^2)/(2h)$ belongs to $K(y)$. It follows
that $[K(y):K]=4$. The monic quartic \eqref{eq:quartic} is its minimal
polynomial, and Gauss's lemma finishes the proof.
\end{proof}

\subsection{Infinitesimal independence}
On a regular branch write \eqref{eq:quartic} as
$w=\phi(u;a,b,R,S)$. In the applications both circle-coordinate
derivatives are nonzero. Indeed, on the selected pair curve,
\[
 F_w=-8h(v-b)(w-a),
\]
so the retained regularity inequalities imply $F_w\ne0$.

\begin{lemma}\label{lem:independence}
For every fixed $a\ne b$ and $R,S>0$, the four analytic functions
$\phi_a,\phi_b,\phi_R,\phi_S$ are linearly independent on every
nonempty regular branch interval. Consequently their Wronskian of
orders zero through three is not identically zero there.
\end{lemma}
\begin{proof}
Suppose that a constant linear combination of the four functions
vanishes. Implicit differentiation gives $\phi_p=-F_p/F_w$.
Therefore the polynomial $Q=\sum_p c_pF_p$ vanishes on a branch
interval of the irreducible quartic. It follows that $F$ divides $Q$.
One way to justify this last assertion is to complexify the regular
analytic branch. The identity theorem makes $Q$ vanish on an open
part of that complex curve. Two coprime plane polynomials have only
finitely many common zeros, so the irreducible polynomial $F$ must
be a factor of $Q$.
But $Q$ has total degree at most three, so $Q=0$.

The homogeneous cubic parts of $F_a,F_b$, in the translated
coordinates $z,y$, are respectively
\[
 4y(z^2-y^2),\qquad -4z(z^2-y^2).
\]
They are independent, forcing $c_a=c_b=0$. If
$A_0=z^2-y^2+R-S-h^2$, then
\[
 F_R=2A_0,\qquad F_S=-2A_0-4h^2.
\]
As $h\ne0$, their constant linear combination is zero only when
$c_R=c_S=0$. This proves independence. The Wronskian criterion
for analytic functions \cite{BostanDumas} gives the last assertion.
Analyticity is relevant here: the corresponding assertion for
arbitrary smooth functions would not suffice.
\end{proof}

\subsection{Uniform sublevel control}
For a sign $\sigma\in\{-1,1\}$, define
\begin{equation}\label{eq:U}
 U_{a,r,\sigma}(v)=a-\sigma\sqrt{r^2-(a-v)^2}.
\end{equation}
The centers and $v$ range over fixed compact intervals, and
$r,s\in[1,2]$. On a retained real chart we require
\begin{equation}\label{eq:regular-chart-margins}
 \sqrt{r^2-(a-v)^2},\quad \sqrt{s^2-(b-v)^2},\quad
 |v-a|,\quad |v-b|\ \geq c\kappa,
 \qquad 0<\kappa<1.
\end{equation}
In particular, both square roots are real and the first-coordinate
map $u=U_{b,s,\sigma_b}(v)$ has nonzero derivative. We also assume
$|a-b|\geq\tau$, where $0<\tau<1$.

We first record an elementary polynomial estimate, including its
uniform dependence on the coefficients. For a polynomial in one
variable, $\|\operatorname{coeff}Q\|_\infty$ denotes the largest
absolute value of its coefficients in the monomial basis.

\begin{lemma}[Polynomial sublevel estimate]\label{lem:polynomial-sublevel}
Let $J$ be a fixed bounded interval of positive length and $d\geq1$.
If $Q$ is a polynomial of degree at most $d$ and
$\|\operatorname{coeff}Q\|_\infty\geq\lambda>0$, then
\begin{equation}\label{eq:polynomial-sublevel}
 |\{v\in J:|Q(v)|\leq\varepsilon\}|
 \leq \min\big\{|J|,C_{d,J}(\varepsilon/\lambda)^{1/d}\big\}
 \qquad(\varepsilon>0).
\end{equation}
The coefficients of $Q$ may be complex.
\end{lemma}
\begin{proof}
Let $E$ be the indicated set and put $m=|E|$. There is nothing to
prove if $m=0$. Otherwise choose $d+1$ points
$v_0,\ldots,v_d\in E$ with mutual distances at least
$\rho=m/(2d+2)$. Such a choice is possible successively: after at
most $d$ choices the excluded intervals have total length at most
$2d\rho<m$. Lagrange interpolation gives
\[
 Q(v)=\sum_{j=0}^d Q(v_j)
       \prod_{\ell\ne j}\frac{v-v_\ell}{v_j-v_\ell}.
\]
The numerator polynomials have coefficients bounded by a constant
depending only on $d,J$, whereas each denominator has absolute
value at least $\rho^d$. Hence
$\lambda\leq C_{d,J}\varepsilon\rho^{-d}
 \leq C'_{d,J}\varepsilon m^{-d}$.
Rearranging and using $m\leq|J|$ proves the assertion.
\end{proof}

To distinguish the two choices of radius parameters, write
$\mathfrak W_{\square}$ for the determinant in \eqref{eq:jet}
with parameter columns $(a,b,R,S)$, where $R=r^2$, $S=s^2$,
and $\mathfrak W$ for the one with columns $(a,b,r,s)$.
The chain rule, with $u$ fixed, gives the exact identity
\begin{equation}\label{eq:radius-wronskian}
 \mathfrak W(u;a,b,r,s)
 =4rs\,\mathfrak W_{\square}(u;a,b,r^2,s^2).
\end{equation}
Thus changing radius parameters introduces a factor bounded above
and below uniformly on $[1,2]^2$.

\begin{lemma}[Uniformity in both radii]\label{lem:sublevel}
There are constants $\alpha>0$ and $C,K<\infty$, independent of
$\kappa,\tau,\delta$, such that for every admissible choice of
$a,b,r,s$ and every choice of branch signs,
\begin{equation}\label{eq:sublevel}
 \left|\left\{v:
  |\mathfrak W(U_{b,s,\sigma_b}(v);a,b,r,s)|\leq\delta
       \right\}\right|
 \leq C(\kappa\tau)^{-K}\delta^\alpha,
 \qquad 0<\delta<1.
\end{equation}
The set on the left is intersected with the retained real charts
in \eqref{eq:regular-chart-margins}.
\end{lemma}
\begin{proof}
We work first with squared radii and put
\[
 q_a(v)=R-(a-v)^2,\qquad q_b(v)=S-(b-v)^2.
\]
We explain separately the algebraic reduction and the uniform
coefficient estimate.

\smallskip\noindent
\emph{Clearing the graph-coordinate denominators.}
Set $u(v)=U_{b,s,\sigma_b}(v)$ and
$w(v)=U_{a,r,\sigma_a}(v)$. Differentiation at fixed parameters is
given by $\partial_u=u_v^{-1}\partial_v$. If $p$ is one of the four
parameters, differentiation of the graph at fixed $u$ gives
\[
 \phi_p=w_p-w_v\frac{u_p}{u_v}.
\]
These identities, and three further $u$-derivatives, express
$\mathfrak W_{\square}$ rationally in
$v,a,b,R,S,\sqrt{q_a},\sqrt{q_b}$. Its denominators consist only
of fixed powers of the square roots and $v-b$. Therefore, for a
sufficiently large fixed integer $M$, the polynomial
\[
 \mathcal D(v;a,b,R,S)=\big(q_a(v)q_b(v)(v-b)\big)^M
\]
clears all denominators. More explicitly,
\begin{equation}\label{eq:cleared-wronskian}
 \mathcal N:=\mathcal D\,
     \mathfrak W_{\square}(u(v);a,b,R,S)
 =A_0+A_1\sqrt{q_a}+A_2\sqrt{q_b}+A_3\sqrt{q_aq_b},
\end{equation}
where the $A_j$ are polynomials in $v,a,b,R,S$ of fixed degree.
Only finitely many derivatives were used, so $M$ and all degree
bounds are independent of the parameters and regularity margins.

\smallskip\noindent
\emph{Taking a polynomial norm.}
Fix $a\ne b$ and $R,S>0$. Each $q_a,q_b$ has two distinct simple
zeros. Their product cannot be a square in $\C(v)$: if they share
one zero, the two unshared zeros still have odd multiplicity; if
they shared both zeros, their common leading coefficient would
force $q_a=q_b$, hence $a=b$. Thus
\[
 \C(v)(\sqrt{q_a},\sqrt{q_b})/\C(v)
\]
is a degree-four extension, whose automorphisms independently
reverse the two square-root signs. The norm
\begin{equation}\label{eq:norm-polynomial}
 P(v;a,b,R,S)
 =\prod_{\epsilon_a,\epsilon_b\in\{-1,1\}}
     \mathcal N_{\epsilon_a,\epsilon_b}(v;a,b,R,S)
\end{equation}
is a polynomial of fixed degree in all five variables. Here
$\mathcal N_{\epsilon_a,\epsilon_b}$ denotes the corresponding
sign conjugate of \eqref{eq:cleared-wronskian}. The invariance of
the product removes the radicals, and the cleared expression
introduces no denominators. Changing the initial branch signs
only permutes the factors in \eqref{eq:norm-polynomial}.

Lemma~\ref{lem:independence} shows that the graph Wronskian is not
identically zero. Its proof applies equally to complex analytic
branches, so this remains true as a statement in the algebraic
function field even if a particular parameter tuple has no common
real chart. Consequently $\mathcal N$ and each of its conjugates
are nonzero field elements. Their norm is therefore not the zero
polynomial in $v$, for every fixed $a\ne b$ and $R,S>0$.

\smallskip\noindent
\emph{Uniformity as the centers approach each other.}
Write $P(v;\mathbf p)=\sum_{j=0}^{d_0}c_j(\mathbf p)v^j$, with
$\mathbf p=(a,b,R,S)$ and a fixed integer $d_0\geq1$.
On the compact parameter box under consideration, let
$G(\mathbf p)=\sum_j|c_j(\mathbf p)|^2$. This is a continuous
semialgebraic function, and the preceding paragraph proves
\[
 G(\mathbf p)=0\quad\Longrightarrow\quad a=b.
\]
The compact semialgebraic \L{}ojasiewicz inequality
\cite[Theorem~2.12]{Coste}, applied to $G$ and $a-b$, yields
$|a-b|^{L_0}\leq C G(\mathbf p)$ for some fixed $L_0$.
After taking square roots and using equivalence of coefficient
norms in this fixed finite-dimensional space, we obtain
\begin{equation}\label{eq:coefficientlower}
 \|\operatorname{coeff}_vP\|_\infty
 \geq c|a-b|^L\geq c\tau^L
\end{equation}
with fixed constants $c,L>0$.

On a real chart, all the sign conjugates in
\eqref{eq:cleared-wronskian} are bounded above by a fixed constant:
they are polynomial expressions in bounded parameters and bounded
real square roots. The same is true of $\mathcal D$. Accordingly,
\eqref{eq:radius-wronskian} implies
\[
 |\mathfrak W(u(v);a,b,r,s)|\leq\delta
 \quad\Longrightarrow\quad |P(v;a,b,R,S)|\leq C\delta.
\]
Apply Lemma~\ref{lem:polynomial-sublevel} on a fixed interval
containing every permitted $v$, with coefficient lower bound
\eqref{eq:coefficientlower}. The removed measure is at most
$C\tau^{-L/d_0}\delta^{1/d_0}$. This proves
\eqref{eq:sublevel}, with $\alpha=1/d_0$ and a sufficiently large
$K$. For the displayed graph determinant the bound is actually
independent of $\kappa$; the weaker form stated above is convenient
when incorporating regular coordinate factors, whose bounds are
polynomial in $\kappa^{-1}$.
\end{proof}

\begin{remark}\label{rem:fiberwise-uniformity}
The quantifiers in Lemma~\ref{lem:sublevel} are essential: the
estimate holds for each fixed pair of radii, with the same
constants. In particular, it bounds the supremum, over the radii,
of the measure of each fiber's bad set. Thus every discarded
integral with a bounded integrand has the same upper bound before
the radius supremum is taken. This permits arbitrary measurable
radius selectors. It does not require estimating the measure of
the union of the bad sets over all radii. A small exceptional set
in the full angle--radius product space alone would not give the
needed bound. The argument uses no derivatives of a radius selector.
\end{remark}

\section{Bounded-input smoothing by a pair-circle argument}\label{sec:bounded}

\begin{proposition}\label{prop:bounded}
There are $\epsilon>0$ and $C<\infty$ such that, for every unit
interval $I$, every $N\geq2$, and arbitrary bounded functions $f,g$,
\begin{equation}\label{eq:bounded}
 \int_I\M(P_Nf,g)(x)\,dx
 \leq C N^{-\epsilon}\norm f_\infty\norm g_\infty.
\end{equation}
\end{proposition}
\begin{proof}
By translation we may fix $I$ in a bounded coordinate box. It is
enough to prove the estimate with $f,g$ supported in a fixed
enlargement of $I$. The cutoff of $g$ is exact, since every physical
argument in \eqref{eq:average} lies in $I+[-2,2]$. Cut off the
\emph{original} first input before applying $P_N$, with a fixed
positive spatial margin. The discarded projection tail is bounded
on $I+[-2,2]$ by $C_LN^{-L}\norm f_\infty$ for every fixed $L$.
This follows directly from the Schwartz kernel of $P_N$.

Normalize the bounded inputs by one and put $F=P_Nf$. For every
fixed finite $p$ we have
\begin{equation}\label{eq:Fbounds}
 \norm F_\infty\leq C,\qquad \norm F_p\leq C_p,
 \qquad\supp\widehat F\subset\{cN\leq|\xi|\leq C_0N\}.
\end{equation}
In particular, $F$ is still exactly frequency localized.

Linearize the supremum by an arbitrary measurable
$r:I\to[1,2]$ and an arbitrary measurable $q$, supported on $I$,
with $|q|\leq1$. We estimate the resulting bilinear form uniformly.
For a finite radius set this reduction is exact: choose, at each
$x$, the first radius attaining the largest modulus, and then choose
the phase of $q(x)$ so that $q(x)\A_{r(x)}(F,g)(x)$ is nonnegative.
These choices are measurable. Uniformity in the finite set permits
passage to the full supremum at the end of the proof.
Remove angular neighborhoods of width $O(\kappa)$ of the coordinate
axes and the directions $\cos\theta+\sin\theta=0$. Their total
contribution to this form is $O(\kappa)$ by boundedness of the
inputs. Use a smooth partition on the remaining finitely many
signed arcs, with derivatives bounded polynomially in
$\kappa^{-1}$.

Fix one such arc. Changing the second input coordinate to $v$
writes its linearization as
\begin{equation}\label{eq:T}
 T_{F,r}g(x)=\int F(U_{x,r(x),\sigma}(v))B(x,r(x),v)g(v)\,dv.
\end{equation}
The density $B$ is smooth in its displayed variables, compactly
supported, and has finite derivative bounds polynomial in
$\kappa^{-1}$. For example, the arclength factor is
$(2\pi|x-U_{x,r,\sigma}(v)|)^{-1}$ times the angular cutoff.

For completeness, set
\[
 J(v)=\int_I q(x)F(U_{x,r(x),\sigma}(v))B(x,r(x),v)\,dx.
\]
Fubini's theorem and Cauchy--Schwarz give
\[
 \left|\int q(x)T_{F,r}g(x)\,dx\right|^2
 =\left|\int g(v)J(v)\,dv\right|^2
 \leq \norm g_2^2\int|J(v)|^2\,dv.
\]
Here $\norm g_2\leq C$ because $g$ is bounded by one and has fixed
compact support. Expanding the square, applying Fubini again, and
using $|q|\leq1$ gives
\begin{equation}\label{eq:TT}
 \left|\int q(x)T_{F,r}g(x)\,dx\right|^2
 \leq C\iint_{I^2}|K_{F,r}(x,y)|\,dx\,dy,
\end{equation}
where, with an immaterial choice of conjugation convention,
\begin{align}\label{eq:kernel}
 K_{F,r}(x,y)=\int &F(U_{x,r(x),\sigma}(v))
                 \overline{F(U_{y,r(y),\sigma}(v))}\notag\\
       &\times B(x,r(x),v)\overline{B(y,r(y),v)}\,dv.
\end{align}
This step removes an arbitrary bounded second input.

Define $H(u,w)=F(u)\overline{F(w)}$. Then
\begin{equation}\label{eq:Hbounds}
 \norm H_{L^p(\R^2)}=\norm F_p^2\leq C_p,
 \qquad\norm H_\infty\leq C,
 \qquad\supp\widehat H\subset\{c_1N\leq|\zeta|\leq C_1N\}.
\end{equation}
The kernel in \eqref{eq:kernel} is a selected value of the linear
two-parameter family
\begin{equation}\label{eq:S}
 S_{r,s}H(x,y)=\int H(U_{x,r,\sigma}(v),U_{y,s,\sigma}(v))
            B(x,r,v)\overline{B(y,s,v)}\,dv.
\end{equation}
Here $(x,y)$ is a genuinely two-dimensional output. We may dominate
the selected value by $\sup_{r,s\in[1,2]}|S_{r,s}H(x,y)|$.
In particular,
\[
 \iint_{I^2}|K_{F,r}(x,y)|\,dx\,dy
 \leq \int_{I^2}\sup_{\rho,\varsigma\in[1,2]}
                  |S_{\rho,\varsigma}H(x,y)|\,dx\,dy.
\]
The two selected radii $r(x)$ and $r(y)$ are simply two admissible
values in this supremum. No derivative of the selector occurs, and
no restriction of an ambient output norm to a graph is involved.

We check the regularity conditions needed to apply
Lemma~\ref{lem:linear} to good pieces of \eqref{eq:S}. Directly,
\begin{equation}\label{eq:derivativesU}
 (U_{x,r,\sigma})_v
 =-\frac{\sigma(x-v)}{\sqrt{r^2-(x-v)^2}},\qquad
 (U_{x,r,\sigma})_x=1-(U_{x,r,\sigma})_v.
\end{equation}
Consequently the output Jacobian of the pair map is
$(1-U_{x,r}'(v))(1-U_{y,s}'(v))$. The removed antidiagonal directions
are exactly its possible zeros on regular arcs. The axis exclusions
ensure regular graph coordinates and a nonzero curve velocity.
All these bounds and coordinate changes cost powers of
$\kappa^{-1}$.
Indeed, when $v=x-r\sin\theta$ and
$U_{x,r,\sigma}(v)=x-r\cos\theta$, one has
\[
 U_v=-\tan\theta,\qquad
 U_x=\frac{\cos\theta+\sin\theta}{\cos\theta}.
\]
The retained angular cutoffs give lower bounds by $c\kappa$ for
$|\sin\theta|$, $|\cos\theta|$, and
$|\cos\theta+\sin\theta|$. Thus both graph velocities and both
factors in the output Jacobian have quantitative lower bounds.
Implicit differentiation supplies all required finite derivative
bounds for the graph-coordinate inverses.

Discard first the output set $|x-y|<\tau$. Its contribution to
the right side of \eqref{eq:TT} is bounded by
$C\kappa^{-C}\tau$. On the complement, the pair curves in
\eqref{eq:S} are the quartics of Section~\ref{sec:geometry} under
the identification $(a,b,R,S)=(y,x,s^2,r^2)$. Permuting the
parameter columns and replacing squared radii by radii introduces
only a sign and factors bounded above and below. Their
corrected cinematic determinant is \eqref{eq:jet}. Insert a smooth
cutoff removing $|\mathfrak W|\leq\delta$ and its transition layer.
By Lemma~\ref{lem:sublevel}, the $v$-measure removed is bounded by
$C(\kappa\tau)^{-C}\delta^\alpha$, uniformly in $x,y,r,s$.
Since $H$ is bounded, the contribution of this part to the kernel
integral obeys the same bound after increasing $C$.

On the good region, all hypotheses of Lemma~\ref{lem:linear} hold
with $A\leq C(\kappa\tau\delta)^{-C}$. A partition into small
regular coordinate boxes has polynomial cardinality and polynomial
finite derivative bounds. Applying the lemma, fixing for instance
$p=100$, and then using H\"older's inequality on $I^2$, we obtain
fixed $c_1>0$ and $C_0,C_2<\infty$ such that
\begin{equation}\label{eq:balance}
 \iint_{I^2}|K_{F,r}|
 \leq C\left[
  \kappa^{-C_0}\tau
  +(\kappa\tau)^{-C_0}\delta^\alpha
  +(\kappa\tau\delta)^{-C_2}N^{-c_1}\right].
\end{equation}
This bound is uniform in the measurable selector $r$.

Here the use of frequency-dependent cutoffs is legitimate for the
following reason. Take smooth functions that vanish for
$|x-y|\leq\tau$ and $|\mathfrak W|\leq\delta$ and equal one when
the respective quantities exceed $2\tau$ and $2\delta$. Every fixed
number of derivatives of the resulting amplitude is bounded by a
power of $(\kappa\tau\delta)^{-1}$. On sufficiently small boxes of
side length $(\kappa\tau\delta)^C$ the graph coordinates are regular
and their determinants retain half their lower bounds. Covering the
five variables $(x,y,r,s,v)$ by these boxes costs only a polynomial
number of pieces. Lemma~\ref{lem:linear} therefore applies with the
same positive frequency exponent on every piece. The transition
layers are included in the preceding exceptional-set bounds.

Choose fixed $L>C_0+1$ and then fixed $J$ such that
$J\alpha>C_0(1+L)+1$. Set $\tau=\kappa^L$ and
$\delta=\kappa^J$. Finally take $\kappa=N^{-b}$, where $b>0$ is
small enough that
\[
 bC_2(1+L+J)<c_1/2,\qquad b\leq c_1/2.
\]
Then \eqref{eq:balance} is $O(N^{-b})$. Taking the square root in
\eqref{eq:TT}, and adding the removed angular contribution, proves
\eqref{eq:bounded}, with $\epsilon=b/2$ after decreasing it if
necessary. The finite arc partition and the original spatial tails
do not alter this positive exponent.
All constants $C_0,C_2,\alpha,c_1,L,J$ are fixed before choosing $b$.
For example, any sufficiently small positive number satisfying
$b<c_1/[2(C_2(1+L+J)+1)]$ works. The finitely bounded range of $N$
for which an angular margin is not yet small is covered by the
elementary bounded-input estimate. No numerical optimality is claimed.
\end{proof}

\section{Interpolation, localization, and translation differences}\label{sec:completion}

We include the remaining steps to keep the passage to general inputs
and to the full radius supremum explicit.

\subsection{An elementary endpoint}
Let
\[
 k_t(u)=\frac{\mathbf 1_{\{|u|<t\}}}{\pi\sqrt{t^2-u^2}}.
\]
Angular Cauchy--Schwarz and H\"older's inequality imply
\begin{equation}\label{eq:elementary}
 |\A_t(F,G)(x)|
 \leq(k_t*|F|^2(x))^{1/2}(k_t*|G|^2(x))^{1/2}
 \leq C\norm F_6\norm G_6,
\end{equation}
uniformly for $1\leq t\leq2$, since the $L^{3/2}$ norms of $k_t$
are uniformly bounded. The multiplier $P_N$ is uniformly bounded
on $L^6$, so \eqref{eq:elementary} is also uniform after projecting
the first input.

For a measurable radius selector, Proposition~\ref{prop:bounded}
is an $L^\infty\times L^\infty\to L^1(I)$ bound of size
$N^{-\epsilon}$, while \eqref{eq:elementary} is an
$L^6\times L^6\to L^\infty(I)$ bound. Bilinear complex
interpolation, with weight $1/4$ on the first endpoint, gives
\begin{equation}\label{eq:local8}
 \norm{\M(P_Nf,g)}_{L^4(I)}
 \leq CN^{-\epsilon/4}\norm f_8\norm g_8.
\end{equation}
Indeed $(3/4)/6=1/8$ and $(1/4)/1=1/4$. The constants are
independent of the selector. Initially, the supremum in
\eqref{eq:local8} may be taken over an arbitrary finite radius set.

\subsection{Spatial summation}
To globalize \eqref{eq:local8}, let $I_k=[k,k+1)$ and use fixed-size
enlargements of these intervals for the original inputs. For a
sufficiently large fixed $A_0$, put
\[
 F_k=\left(\int(1+|u-k|)^{-A_0}|f(u)|^8\,du\right)^{1/8},
 \qquad
 G_k=\left(\int(1+|u-k|)^{-A_0}|g(u)|^8\,du\right)^{1/8}.
\]
Then $\sum_kF_k^8\leq C\norm f_8^8$ and similarly for $G_k$.
The near-input part of \eqref{eq:local8} is bounded by
$CN^{-\epsilon/4}F_kG_k$. The far part of the original first input,
after projection, has pointwise size at most $C_MN^{-M}F_k$ on
$I_k+[-2,2]$ for every fixed $M$. This is the Schwartz-kernel tail
bound followed by H\"older's inequality. The other factor is
estimated by the local version of \eqref{eq:elementary}; its
local $L^6$ norm is bounded by $CG_k$. Thus
\[
 \norm{\M(P_Nf,g)}_{L^4(I_k)}
 \leq CN^{-\epsilon/4}F_kG_k.
\]
Summing fourth powers and using Cauchy--Schwarz in $k$ proves
Proposition~\ref{prop:frequency}, with $s=\epsilon/4$.
The same localization argument applied directly to
\eqref{eq:elementary} gives
\begin{equation}\label{eq:ordinary}
 \norm{\M(f,g)}_4\leq C\norm f_8\norm g_8.
\end{equation}

\subsection{Littlewood--Paley summation}
Take an inhomogeneous smooth Littlewood--Paley decomposition
$I=P_{\leq1}+\sum_{j\geq1}P_{2^j}$. We may decrease $s$ so that
$0<s<1$. If $\widetilde P_N$ is an enlarged annular cutoff equal
to one on the support of $P_N$, the difference-kernel estimate gives
\begin{equation}\label{eq:differencekernel}
 \norm{P_N(I-\tau_h)f}_8
 \leq C\min(1,N|h|)\norm f_8.
\end{equation}
Indeed, if $K_N(x)=N\check p(Nx)$ is the projection kernel, then
\[
 \norm{K_N-\tau_hK_N}_1
 \leq\min\bigl(2\norm{K_N}_1,
                  |h|\norm{K_N'}_1\bigr)
 \leq C\min(1,N|h|).
\]
Young's inequality proves \eqref{eq:differencekernel}; the same
argument gives $C|h|$ for the low-pass multiplier.
Apply \eqref{eq:globalfrequency} with $\widetilde P_N$ to the first
input $P_N(I-\tau_h)f$. For $0<|h|\leq1$, Minkowski's inequality
and \eqref{eq:differencekernel} yield
\begin{align*}
 \norm{\M((I-\tau_h)f,g)}_4
 &\leq C\left(|h|+\sum_{j\geq1}2^{-js}
                       \min(1,2^j|h|)\right)\norm f_8\norm g_8\\
 &\leq C_s|h|^s\norm f_8\norm g_8.
\end{align*}
To verify the last sum directly, choose $J\geq0$ with
$2^J\leq|h|^{-1}<2^{J+1}$. The part $j\leq J$ is at most
$C_s|h|2^{J(1-s)}\leq C_s|h|^s$, and the part $j>J$ is at most
$C_s2^{-Js}\leq C_s|h|^s$.
The low-frequency term uses \eqref{eq:ordinary} and the low-pass
analogue of \eqref{eq:differencekernel}. Symmetry of the circle
gives the estimate for a difference in the second input. For
$|h|\geq1$, \eqref{eq:ordinary} suffices. This proves
Theorem~\ref{thm:main} with $\eta=s$.

Finally, the frequency series is initially identified with the
operator for smooth compactly supported inputs. Its output series
is absolutely convergent in $L^4$ by
Proposition~\ref{prop:frequency}. Bound \eqref{eq:ordinary} and
density identify the limit for all $L^8$ inputs. On a fixed compact
output interval, $L^8$ convergence of localized inputs implies
$L^6$ convergence, so \eqref{eq:elementary} gives uniform convergence
of the averages in $(x,t)$. Approximation by smooth inputs therefore
shows that the averages have a jointly continuous representative
there. The supremum over a countable dense subset of $[1,2]$ agrees
with the full supremum, completing the measurability and
linearization justification.

\subsection{Propagation to other exponent triples}
We state explicitly why a single triple answers the existence
question. Write
\[
 \mathcal R=\left\{(1/p,1/q,1/r):
  1\leq p,q\leq\infty,\ 0<r\leq\infty,\
  \norm{\M(f,g)}_r\leq C_{p,q,r}\norm f_p\norm g_q\right\},
\]
with the usual quasi-norm interpretation if $r<1$. This is a convex
set by multilinear interpolation. We use strong boundedness here;
no weak-type endpoint is included without an additional argument.

\begin{proposition}\label{prop:interpolation}
Let $a_0=(1/8,1/8,1/4)$ and let $a_1=(1/p_1,1/q_1,1/r_1)$ belong
to $\mathcal R$. For $0<\theta<1$ define
\begin{equation}\label{eq:interpolatedtriple}
 (1/p_\theta,1/q_\theta,1/r_\theta)
     =\theta a_0+(1-\theta)a_1.
\end{equation}
Then the two translation-continuity estimates hold at this triple
with power $\theta\eta$, where $\eta$ is as in
Theorem~\ref{thm:main}. Consequently they hold at every point in
the relative interior of $\mathcal R$, with a power depending on
the point.
\end{proposition}
\begin{proof}
Fix $0<|h|\leq1$ and first take a finite radius set. For any
measurable selector $\rho(x)$ from that set, the operator
\[
 B_{\rho,h}(f,g)(x)=\A_{\rho(x)}((I-\tau_h)f,g)(x)
\]
is bilinear. Its norm at $a_0$ is at most $C|h|^\eta$, uniformly
in $\rho$, by Theorem~\ref{thm:main}. At $a_1$, the triangle
inequality in the first input and translation invariance of its
norm give the bound $2C_{p_1,q_1,r_1}$. Multilinear complex
interpolation, in the form of \cite[Proposition~2.2]{BOS}, gives
\[
 \norm{B_{\rho,h}(f,g)}_{r_\theta}
 \leq C_\theta|h|^{\theta\eta}
                \norm f_{p_\theta}\norm g_{q_\theta}.
\]
That interpolation theorem also permits quasi-Banach target
Lebesgue spaces. The bound is independent of the radius set and
of $\rho$, so measurable selection gives the corresponding finite
maximal estimate. Initially take smooth compactly supported inputs.
For them the averages are continuous in the radius, and increasing
finite subsets of a dense radius set recover the full supremum.
Monotone convergence gives the desired estimate. Ordinary maximal
boundedness at the interpolated triple then extends it by density
to all inputs. The input exponents at that triple are finite because
$\theta>0$. The second-input difference is handled in exactly the
same way.

If $a$ is in the relative interior of $\mathcal R$, then for all
sufficiently small $\theta>0$ the point
\[
 a_1=\frac{a-\theta a_0}{1-\theta}
     =a+\frac{\theta}{1-\theta}(a-a_0)
\]
belongs to $\mathcal R$. Indeed $a-a_0$ is a direction in its
affine hull, and a relative neighborhood of $a$ is contained in
$\mathcal R$. Apply \eqref{eq:interpolatedtriple} to this $a_1$.
\end{proof}

\begin{proof}[Proof of Corollary~\ref{cor:holdercontinuity}]
The global bilinear circular maximal bounds of
\cite[Theorem~1.3(ii)]{CZ} and \cite[Theorem~1]{DR} imply
\[
 \norm{\M(f,g)}_{pq/(p+q)}\leq C_{p,q}\norm f_p\norm g_q
 \qquad(2<p,q<\infty).
\]
The local operator is pointwise dominated by the global positive
maximal operator considered in those results. For fixed finite
$p,q>2$, choose $\theta>0$ so small that
\[
 \frac1{p_1}=\frac{1/p-\theta/8}{1-\theta},\qquad
 \frac1{q_1}=\frac{1/q-\theta/8}{1-\theta}
\]
both belong to $(0,1/2)$. Define $1/r_1=1/p_1+1/q_1$.
These exponents satisfy \eqref{eq:interpolatedtriple} with
$(p_\theta,q_\theta,r_\theta)=(p,q,pq/(p+q))$.
Proposition~\ref{prop:interpolation} proves the corollary with
$\eta_{p,q}=\theta\eta$.
\end{proof}

\appendix
\section{Quantitative dependence in the linear smoothing input}\label{app:quantitative}

This appendix supplies the quantitative version of the linear result
used in Lemma~\ref{lem:linear}. The dependence on $A$ in
\eqref{eq:lineardecay} is obtained by following the normalization in
the proof of \cite{CGY}; it is not stated separately in its
Theorem~1.5. We explain both the normalization and the frequency
exponents, since qualitative dependence on the curve would not
suffice for the frequency-dependent cutoffs in
Section~\ref{sec:bounded}.

Fix $p>12$ and $\varepsilon_*>0$. All numbers of derivatives and all
error orders in what follows are chosen in terms of these two
numbers before $A$ or $N$ is chosen. Constants $K$ below are finite
and may increase. We shall first prove the joint estimate
\begin{equation}\label{eq:joint}
 \norm{\mathcal S_vH(z)}_{L^p_{z,v}}
 \leq C_{p,\varepsilon_*} A^K
             N^{-3/p+\varepsilon_*}\norm H_p.
\end{equation}
The corresponding derivative estimates will then give the maximal
bound by Sobolev embedding in the two parameters.

\subsection{Quantitative coordinate changes}
We use the following elementary observations repeatedly. If a
matrix of fixed size has entries bounded by $A^C$ and determinant
of magnitude at least $A^{-C}$, its inverse has entries bounded by
$A^{C'}$, by the cofactor formula. Differentiating the identity
$B B^{-1}=I$ gives the same statement for any prescribed finite
number of derivatives. The quantitative inverse function theorem
therefore gives inverse coordinate maps, with polynomial derivative
bounds, on balls of radius $A^{-C'}$. A fixed compact box can be
covered by a polynomial number of such balls. Partitions of unity
on them also have polynomial finite derivative bounds.

Apply these observations to the graph normalizations in
Lemma~\ref{lem:linear}. On each resulting patch write
\[
 \Gamma(X,Y;\theta;v)=(X-\theta,h(X,Y;\theta;v)),
 \qquad Z=(X,Y,v_1,v_2).
\]
The corrected determinant \eqref{eq:correctedH3} retains a lower
bound by an inverse power of $A$. To avoid a potentially small
pure curve derivative, perform the fixed input-coordinate shear
\begin{equation}\label{eq:quant-shear}
 T(u,w)=(u,w-c u^2),\qquad c=A_1+1,
\end{equation}
where $A_1\leq A^C$ bounds $|h_{\theta\theta}|$ on the patch.
The transformed graph is $h-c(X-\theta)^2$, so its second
$\theta$ derivative is at most $-A_1-2$. Every $\theta$ derivative
of the added term is annihilated by
$\partial_\theta+\partial_X$; its $Y,v_1,v_2$ derivatives vanish.
Thus \eqref{eq:correctedH3} is unchanged, exactly. Both $T$ and
$T^{-1}$ have polynomial finite derivative bounds and unit Jacobian.
After this shear we may use the second-derivative stationary-phase
case throughout, without selecting a smaller uncontrolled derivative.

\subsection{Frequency localization under the input shear}
The shear is a change of physical coordinates, so it does not
preserve exact Fourier support. We record the required replacement.
Let $T$ be any of the localized input diffeomorphisms just described,
and let $\chi$ be a compactly supported cutoff on its coordinate
patch. For an annular input of frequency $N$, the Fourier kernel of
$H\mapsto\chi(H\circ T)$ is
\begin{equation}\label{eq:composition-kernel}
 I(\xi,\zeta)=\int\chi(u)
          e^{i[T(u)\cdot\zeta-u\cdot\xi]}\,du,
 \qquad |\zeta|\simeq N.
\end{equation}
The singular values of $DT$ lie between powers $A^{-C}$ and $A^C$.
Hence, for a sufficiently enlarged frequency window
\[
       N A^{-C_1}\leq|\xi|\leq N A^{C_1},
\]
the phase in \eqref{eq:composition-kernel} has no stationary point
outside the window. Its gradient there is bounded below by
$A^{-C}(N+|\xi|)$. Repeated integration by parts gives
\begin{equation}\label{eq:composition-offband}
 |\partial_\xi^\mu\partial_\zeta^\nu I(\xi,\zeta)|
 \leq C_{L,\mu,\nu} A^K(N+|\xi|)^{-L}
\end{equation}
for any fixed orders $L,\mu,\nu$, provided the initially chosen
number of derivatives is large enough.

Insert an input annular cutoff and an output dyadic cutoff of
radius $\rho\geq1$. Integration in the frequency variables converts
\eqref{eq:composition-offband} into a physical kernel whose two
Schur norms are bounded by
\begin{equation}\label{eq:composition-schur}
       C A^K N^2\rho^2(N+\rho)^{-L}.
\end{equation}
For the low-frequency output block, replace $\rho$ by $1$.
Consequently the sum of the off-window pieces has $L^p$ norm at
most $C_M A^K N^{-M}\norm H_p$, for any prescribed fixed $M$.
We may also multiply \eqref{eq:composition-schur} by the Bernstein
factor $\rho^{2/p}$ and sum. This bounds that same error in
$L^\infty$ by $C_M A^K N^{-M}\norm H_p$. Its contribution to a
curve average, even after the parameter supremum, is therefore
negligible by the amplitude's $L^1$ bound. No unquantified maximal
operator estimate is needed to dispose of the error.

The retained window consists of $O(1+\log A)$ dyadic annuli.
The $L^p$ norm of composition is bounded by a power of $A$, and
each annular projection is uniformly $L^p$ bounded. It is therefore
enough to establish the desired estimates at each retained frequency,
which differs from $N$ by at most a power of $A$. The finite sum
changes only the polynomial prefactor. This also handles the
off-annular tail introduced by the compact input cutoff.

\subsection{A normalized finite-derivative class}
For completeness, we spell out the normalization relevant to the
external decoupling theorem. Discard nonstationary frequency sectors
by integration by parts and put $\omega=\xi/\eta$ on a remaining
sector. The stationary point $\theta_*(Z,\omega)$ is determined by
\[
 h_\theta(Z,\theta_*)=\omega.
\]
The lower bound for $|h_{\theta\theta}|$ and the preceding inverse
function bounds control its finite derivatives polynomially. The
stationary phase is
\[
 \Psi(Z,\omega)=(X-\theta_*(Z,\omega))\omega
                        +h(Z,\theta_*(Z,\omega)).
\]
For a fixed $\omega_0$, introduce the four joint-output coordinates
\begin{equation}\label{eq:phase-jet-map}
 Z\longmapsto\left(
 \Psi(Z,\omega_0),\partial_\omega\Psi(Z,\omega_0),
 \partial_\omega^2\Psi(Z,\omega_0),
 \partial_\omega^3\Psi(Z,\omega_0)\right).
\end{equation}
The determinant computation in \cite[Section~8.2, (8.22)--(8.34)]{CGY}
identifies this Jacobian, up to sign, with the corrected cinematic
determinant times a power of $\partial_\omega\theta_*$.
Both it and its inverse are consequently controlled by powers of
$A$. These are changes in all four joint-output variables, which
are permissible for the joint norm in \eqref{eq:joint}.

Normalize the matrix of \eqref{eq:phase-jet-map}, restrict the
joint-output coordinates to sufficiently small boxes, and restrict
$\omega$ to intervals of length $A^{-C}$. Taylor expansion through
degree three, followed by the corresponding anisotropic dilations
of the four jet coordinates, puts the phase into a fixed bounded
finite-derivative class with moment-cone jet
$(1,t,t^2,t^3)$. The remainder and all its required derivatives are
made smaller than a fixed normalization constant by taking $C$
large enough. The number of pieces, coordinate Jacobians, amplitude
derivatives, and changes of frequency scale are all bounded by
fixed powers of $A$.

The needed uniformity of decoupling over this class is the one used
to close its frequency induction. For quadratic cones it is explicit
in the type-$\mathbf 1$ formulation of
\cite[Sections~2.3--2.7, especially Lemma~2.3 and Remark~2.4]{BHS}.
The degree-three version used in \cite[Lemma~7.2]{CGY} has the same
finite-derivative normalization and bootstrapping argument
\cite{PS}, with moment-curve decoupling \cite{BDG} as its
constant-coefficient input. The derivative and error orders in
these arguments are fixed in terms of $p,\varepsilon_*$.
Thus the normalized induction constant is independent of $A$.
The conditioning just described is performed once, before that
induction; it produces a factor $A^K$, not a new factor of $A$
at each dyadic frequency scale.

\subsection{The two decoupling scales}
We give the frequency calculation explicitly for $d=3$. Write
$\lambda$ for a frequency in the normalized coordinates; thus
$N A^{-C}\leq\lambda\leq N A^C$. Stationary phase in $\theta$
gives an amplitude of order $\lambda^{-1/2}$ with uniform symbol
bounds, up to the polynomial prefactors already recorded. The
degree-three decoupling step down to ratio intervals of length
$\lambda^{-1/3}$ contributes
\[
 \lambda^{\frac13(1-7/p)+\varepsilon_*/2}.
\]
On each such interval the degree-two step, down to length
$\lambda^{-1/2}$, contributes
\[
 \lambda^{(\frac12-\frac13)(1-4/p)+\varepsilon_*/2}.
\]
These are \cite[(7.27)--(7.30)]{CGY}, with $J_3=J_2=\log_2\lambda$;
the nonvanishing second derivative places us directly in their
second-derivative case. The degree-two condition is uniform in the
normalized moment-cone class, as in their (7.28).
More explicitly, on a ratio interval of length
$\Delta=\lambda^{-1/3}$ the cubic term has size
$\lambda\Delta^3=1$ and can be included in a uniformly smooth
amplitude. The quadratic part has effective scale
$\lambda\Delta^2=\lambda^{1/3}$, with a nonzero quadratic
parameter derivative in the normalized jet coordinates. Its
degree-two decoupling factor is therefore
$(\lambda^{1/3})^{(1-4/p)/2}$, which is the displayed second factor.

At the final angular scale, the Taylor remainder of the phase
after its constant and linear angular terms has bounded size and
bounded rescaled derivatives. The resulting cap operator has an
integrable kernel with uniform Schur bounds in appropriate local
output coordinates. This is the final Young-inequality step
\cite[(7.31)]{CGY}. In the variable-coefficient formulation the
first two phase-jet coordinates give the two corresponding input
coordinates; the other two joint-output coordinates are integrated
over bounded sets. Their coordinate changes have already been
included in the polynomial bounds.

Here is an explicit kernel calculation for this last assertion.
At the center $\omega_\Theta$ of a final cap, use
\[
 q_1=\Psi_\omega(Z,\omega_\Theta),\qquad
 q_0=\Psi(Z,\omega_\Theta)-\omega_\Theta
                                      \Psi_\omega(Z,\omega_\Theta)
\]
as two coordinates, completed by the remaining two jet coordinates.
Then
$\eta\Psi=q_1\xi+q_0\eta+\eta R$, where
$R=O((\omega-\omega_\Theta)^2)$. Put
$\Delta=\lambda^{-1/2}$. On the cap $\eta R$ has bounded
derivatives in the rescaled frequency coordinates. Repeated
integration by parts bounds its kernel, with the stationary-phase
factor already extracted, by
\[
 C_M\lambda^2\Delta\left(1+
 \lambda\Delta|q_1-u_1|+
 \lambda|q_0-u_2+\omega_\Theta(q_1-u_1)|\right)^{-M}.
\]
For $M>2$, integration in either $(u_1,u_2)$ or $(q_1,q_0)$
is uniformly bounded. Integration in the two remaining bounded
coordinates, and their controlled Jacobians, gives the asserted
joint-output estimate for each cap.

The final input projections satisfy
\[
 \left(\sum_\Theta\norm{H_\Theta}_p^p\right)^{1/p}
 \leq C_p\norm H_p\qquad(p\geq2).
\]
Indeed, Plancherel gives the assertion at $p=2$ by bounded overlap,
whereas the rescaled multiplier kernels give the
$L^\infty\to\ell^\infty(L^\infty)$ bound; interpolation gives
the displayed inequality. Any frequency-only phase introduced in
an intermediate normalization is harmless for decoupling itself,
but is retained in the final cap kernel calculation. We do not
regard a general unimodular Fourier multiplier as an $L^p$ isometry.

Combining the two scales with stationary phase gives
\begin{align*}
 &\lambda^{-1/2}
 \lambda^{\frac13(1-7/p)}
 \lambda^{\frac16(1-4/p)}
 \lambda^{\varepsilon_*}
   =\lambda^{-3/p+\varepsilon_*}.
\end{align*}
This proves \eqref{eq:joint}, after returning to the original
coordinates and absorbing the changes of scale into $A^K$.
Finite stationary-phase expansions and Fourier-series expansions
of amplitudes use only the predetermined finite derivatives;
their remainder estimates have the same polynomial dependence.

\subsection{From the joint norm to the maximal norm}
The same proof applies after a fixed integer number $j$ of
derivatives in $v$: differentiating the original oscillatory
formula multiplies its order by at most $A^K N^j$. Therefore
\[
 \norm{\partial_v^\beta\mathcal S_vH}_{L^p_{z,v}}
 \leq C A^K N^{|\beta|-3/p+\varepsilon_*}\norm H_p.
\]
Interpolation between integer Sobolev orders gives the analogous
estimate of fractional order $s$. Apply
$W^{s,p}(\R^2_v)\hookrightarrow L^\infty(\R^2_v)$ with
$s=2/p+\varepsilon_*$, using a compact parameter cutoff. Taking
the $L^p_z$ norm yields
\[
 \norm{\sup_v|\mathcal S_vH|}_p
 \leq C A^K N^{-1/p+2\varepsilon_*}\norm H_p.
\]
Choose $\varepsilon_*=1/(4p)$ to obtain
\eqref{eq:lineardecay}. The argument applies when the normalized
frequency is at least one. In the remaining range $N\leq A^C$,
the elementary bandlimited estimate
$\norm H_\infty\leq C N^{2/p}\norm H_p$, followed by the
amplitude's $L^1$ bound on a compact output box, proves the same
claim after increasing $K$.

All cutoffs, coordinate patches, and error estimates in this
argument have polynomial cost in $A$. The explicit
shrinking-curvature bound in \cite[(3.19)]{CGY} is consistent with
this bookkeeping. The preceding argument establishes precisely
the uniformity needed when $A$ itself is a small power of $N$.

\bigskip
\noindent School of Mathematical Sciences, Beihang University, Beijing China

\noindent E-mail adress: \texttt{weixuyi@buaa.edu.cn}


\begin{thebibliography}{99}

\bibitem{BOS}
J.-G. Bak, D. M. Oberlin, and A. Seeger,
\emph{Restriction of Fourier transforms to curves and related
oscillatory integrals},
Amer. J. Math. \textbf{131} (2009), no.~2, 277--311.
\href{https://arxiv.org/abs/math/0612752}{arXiv:math/0612752}.

\bibitem{BHS}
D. Beltran, J. Hickman, and C. D. Sogge,
\emph{Variable coefficient Wolff-type inequalities and sharp local
smoothing estimates for wave equations on manifolds},
Anal. PDE \textbf{13} (2020), no.~2, 403--433.
\href{https://arxiv.org/abs/1801.06910}{arXiv:1801.06910}.

\bibitem{Borges}
T. Borges,
\emph{Survey on bilinear spherical averages and associated maximal
operators},
J. Geom. Anal. \textbf{36} (2026), article~94.
\href{https://doi.org/10.1007/s12220-026-02336-w}
{doi:10.1007/s12220-026-02336-w}.

\bibitem{BFOPZ}
T. Borges, B. Foster, Y. Ou, J. Pipher, and Z. Zhou,
\emph{Sparse bounds for the bilinear spherical maximal function},
J. Lond. Math. Soc. (2) \textbf{107} (2023), no.~4, 1409--1449.
\href{https://doi.org/10.1112/jlms.12715}{doi:10.1112/jlms.12715}.

\bibitem{BostanDumas}
A. Bostan and P. Dumas,
\emph{Wronskians and linear independence},
Amer. Math. Monthly \textbf{117} (2010), no.~8, 722--727.
\href{https://arxiv.org/abs/1301.6598}{arXiv:1301.6598}.

\bibitem{Bourgain}
J. Bourgain,
\emph{Averages in the plane over convex curves and maximal operators},
J. Analyse Math. \textbf{47} (1986), 69--85.
\href{https://doi.org/10.1007/BF02792533}{doi:10.1007/BF02792533}.

\bibitem{BDG}
J. Bourgain, C. Demeter, and L. Guth,
\emph{Proof of the main conjecture in Vinogradov's mean value theorem
for degrees higher than three},
Ann. of Math. (2) \textbf{184} (2016), no.~2, 633--682.
\href{https://doi.org/10.4007/annals.2016.184.2.7}
{doi:10.4007/annals.2016.184.2.7}.

\bibitem{CGY}
M. Chen, S. Guo, and T. Yang,
\emph{A multi-parameter cinematic curvature},
\href{https://arxiv.org/abs/2306.01606v7}{arXiv:2306.01606v7},
27 January 2025. The corrected hypothesis (H3) in this version is used.

\bibitem{CZ}
M. Christ and Z. Zhou,
\emph{A class of singular bilinear maximal functions},
J. Funct. Anal. \textbf{287} (2024), no.~8, article~110572.
\href{https://doi.org/10.1016/j.jfa.2024.110572}
{doi:10.1016/j.jfa.2024.110572}.

\bibitem{Coste}
M. Coste,
\emph{An introduction to semialgebraic geometry},
lecture notes, 2002.
\href{https://perso.univ-rennes1.fr/michel.coste/polyens/SAG.pdf}
{Available from the author's webpage}.

\bibitem{DR}
G. Dosidis and J. P. G. Ramos,
\emph{The multilinear spherical maximal function in one dimension},
Proc. Edinb. Math. Soc. (2) \textbf{67} (2024), no.~4, 1045--1059.
\href{https://doi.org/10.1017/S0013091524000191}
{doi:10.1017/S0013091524000191}.

\bibitem{MSS}
G. Mockenhaupt, A. Seeger, and C. D. Sogge,
\emph{Local smoothing of Fourier integral operators and
Carleson--Sj\"olin estimates},
J. Amer. Math. Soc. \textbf{6} (1993), no.~1, 65--130.
\href{https://doi.org/10.1090/S0894-0347-1993-1168960-6}
{doi:10.1090/S0894-0347-1993-1168960-6}.

\bibitem{PS}
M. Pramanik and A. Seeger,
\emph{$L^p$ regularity of averages over curves and bounds for
associated maximal operators},
Amer. J. Math. \textbf{129} (2007), no.~1, 61--103.
\href{https://arxiv.org/abs/math/0501165}{arXiv:math/0501165}.

\bibitem{Sogge}
C. D. Sogge,
\emph{Propagation of singularities and maximal functions in the plane},
Invent. Math. \textbf{104} (1991), 349--376.
\href{https://doi.org/10.1007/BF01245080}{doi:10.1007/BF01245080}.

\end{thebibliography}
\end{document}